\documentclass[10pt]{article}

\usepackage{amssymb,amsmath,latexsym,amsfonts,graphicx,hyperref}      
                 
   \usepackage{epsfig} 
\usepackage{amssymb,amsmath,amsfonts,epsfig,color,placeins,subfigure,xcolor,amsthm}
\usepackage{enumerate}

\usepackage{amssymb,amsmath,latexsym,color,dsfont,amsfonts}
\usepackage{graphicx}

\definecolor{MyDarkRed}{rgb}{0.8,0.08,0}
\definecolor{HisDarkblack}{rgb}{0,0.08,0.8}

\newtheorem{prop}{Proposition}
\newtheorem{lemma}{Lemma}
\newtheorem{ass}{Assumption}
\newtheorem{theo}{Theorem}
\newtheorem{corol}{Corollary}
\newtheorem{remark}{Remark}
\newtheorem{deff}{Definition}

\begin{document}

\title{\LARGE 
Backstepping control for the sterile mosquitoes technique: stabilization of extinction equilibrium}

\author{Andrea Cristofaro\thanks{Dipartimento di Ingegneria Informatica, Automatica e Gestionale ``A. Ruberti", Sapienza University of Rome, Italy. \texttt{andrea.cristofaro@uniroma1.it}}$\ $ and Luca Rossi\thanks{Dipartimento di Matematica ``G. Castelnuovo", Sapienza University of Rome, Italy. \texttt{l.rossi@uniroma1.it}}}



\pagestyle{empty} 

\maketitle

\thispagestyle{empty}

\begin{abstract}            
The control of a mosquito population using the sterile insect technique is considered. Building on a model-based approach, where the control input is the release rate of sterilized males, we propose a non-negative backstepping control law capable of globally stabilizing the extinction equilibrium of the system. A simulation study supports and validates the theoretical findings, showing the efficacy of the approach both on a reduced model, used for control design, and on a complete model of the mosquito population dynamics.

\end{abstract}

\section{Introduction}
 Among several approaches against mosquito-borne diseases, such as malaria and dengue, sterile insect technique (SIT) is a promising and fascinating one. The basic idea is to release  sterilized male mosquitoes, \textcolor{black}{typically after irradiation}, so that a substantial portion of the laid eggs by females are not fertile 
 {\cite{knipling1968genetic,harris2011field,dyck2021sterile}.} While SIT has been successfully used against different insect pests, it is particularly well suited for the case of mosquitoes because the males do not feed on blood and are completely harmless for humans.
Nevertheless, the performance and the efficacy of the method are quite complex to assess, due to several exogenous variables involved, such as local climate and environment, as well as dependency on other factors such as rate and amount of sterile males being released and spatial diffusion phenomena \cite{almeida2023analysis,leculier2023control}. 

Bearing this in mind, the goal of the present note is then to further develop a control-theoretic framework to analyze and implement effective strategies for the sterile insect technique applied to a mosquito population. \textcolor{black}{Impulsive feedback control has been proposed in \cite{bliman2019implementation}, while optimal control strategies have been considered in \cite{almeida2019mosquito,almeida2022optimal}. We propose here a backstepping control design setup based on a reduced model. Related results can be found also in the recent paper \cite{bidi2023global}.} The form of the controller, due to the particular structure of the system, is non trivial, and additional efforts are needed to show that the proposed feedback is always non-negative, thus making it consistent with the biological interpretation.
Indeed, the backstepping control law is proven to globally exponentially stabilize the extinction equilibrium for the reduced system, and based on the simulation results, such property seems to be retained when the same controller is applied to the complete population model. Simulations also highlight an interesting point concerning the robustness of the approach: in fact, the backstepping control law appears to have an inherent ability of being robust to uncertainties in the model parameters. \smallskip\\
The paper is structured as follows: the model is presented in Section \ref{sec:model}, while the main contributions are reported in Section \ref{sec:control}, which covers in particular the design of the proposed backstepping control law and the analysis of the system stability. Simulation results are proposed in Section~\ref{sec:simul} and, finally, conclusions are drawn in Section~\ref{sec:conclusion}.

\section{Mosquito population modeling}\label{sec:model}
Consider a mosquito population divided in four compartments, representing {densities [\#]} of: individuals in aquatic phase (larvae) $E$, adult males $M$, fertilized adult females $F$ and sterilized adult males $M_s$.  \textcolor{black}{We gather from \cite{strugarek2019use} and  \cite[Section 2]{almeida2022optimal} the key
ingredients required presently. The population dynamics is then governed by the following compartmental model  \cite{haddad2010nonnegative}:}
\begin{equation}\label{eq:fullmodel}
\begin{array}{l}
\dot{E}=\beta_E F\left(1-\frac{E}{k}\right)-(\nu_E+\delta_E)E\\
\dot{M}=(1-\nu)\nu_EE-\delta_MM\\
\dot{F}=\nu \nu_EE\frac{M}{M+\gamma_s M_s}-\delta_FF\\
\dot{M_s}=u-\delta_s M_s
\end{array}
\end{equation}
where
\begin{itemize}
\item $\beta_E>0$ is the oviposition rate {[\#/day]}
\item $\delta_*>0$ indicates the death rate for individuals belonging to class $*\in\{E,M,F,s\}$ {[\#/day]}
\item $\nu_E$ is the hatching rate for eggs {[\#/day]}
\item $\nu\in(0,1)$ is the probability that a pupa\footnote{{Pupa refers to the intermediate stage between a larva and a mature insect.}} becomes a female mosquito {[adimensional]}
\item $k>0$ is the 
environmental capacity for eggs\footnote{It refers to the maximum density of eggs.} {[\#]}
\item $\gamma_s>0$ is the index of preference of females for sterile males against fertile males\footnote{In particular, the probability for a female to mate with a fertile male is $M/(M+\gamma_sM_s)$} {[adimensional]}
\item $u$ is the control function, representing the rate of delivery of sterile males {[\#/day]}
\end{itemize}
The following {modelling} assumption is made, related to the fact that sterile males are somehow less fit than both females and fertile males, and that females generate more than one female in their lifespan. 
\begin{ass}\label{ass:offspring}
The model parameters satisfy
$$
{\delta_s>\max\{\delta_F,\delta_M\}},\quad R_0:=\frac{\nu\beta_E\nu_E}{\delta_F(\nu_E+\delta_E)}>1.
$$
$R_0$ is called basic offspring number.
\end{ass}
\textcolor{black}{The first condition in Assumption~\ref{ass:offspring} is motivated by
empirical observations, see e.g.~\cite{Ernawan22}. Moreover, it addresses the most challenging case 
--where sterilized males are weaker than other mosquitoes-- 
since the greater the number of sterilized males,
the higher the probability of extinction of the population.}
{The second condition in Assumption~\ref{ass:offspring} implies that the population is naturally bounded to persist, whereas $R_0<1$ would imply natural extinction.}

\textcolor{black}{The typical lifespans of eggs (3-4 days) and males (9-12 days) 
may be considerably shorter than the one of females (up to 6 weeks, depending on the species and the environmental conditions) 
\cite{lifespan, sowilem2013life}.
This motivates the following simplification of the model: introducing the time rescaling $\tilde{t}=t/\varepsilon$ and the new functions\footnote{Note that, by definition, $\tilde{E}\!\left(\tilde{t}\right)=E(t)$ and $\tilde{M}\!\left(\tilde{t}\right)=M(t)$.} $\tilde{E}(t):=E(\varepsilon t)$ and
$\tilde{M}(t):=M(\varepsilon t)$, one derives from 
\eqref{eq:fullmodel} 
$$\frac{d}{d\tilde{t}}\,\tilde{E}\!\left(\tilde t\right)=\varepsilon \frac{d}{dt} E(t),\quad \frac{d}{d\tilde{t}}\,\tilde{M}\!\left(\tilde t\right)=\varepsilon \frac{d}{dt} M(t)$$
whence, letting $\varepsilon\to0$, one formally finds that the compartments $E$ and $M$ are at equilibrium (see also \cite[Section 2]{almeida2019mosquito}):
\begin{equation}\label{eq:comp_eq}
E=\frac{\beta_E F}{\alpha(F)},\quad M=\frac{(1-\nu)\nu_E}{\delta_M} E,
\end{equation}
with $\alpha(F):=\frac{\beta_E F}{k}+\nu_E+\delta_E$.
Plugging the expressions \eqref{eq:comp_eq} into \eqref{eq:fullmodel} leads}
to the reduced model \cite{almeida2019mosquito}, \cite{almeida2022optimal}
\begin{equation}\label{eq:reduced}
\begin{array}{l}
\dot{F}=g(F,M_s)-\delta_FF\\
\dot{M_s}=u-\delta_sM_s
\end{array}
\end{equation}
where the nonlinear function $g:\mathbb{R}\times\mathbb{R}\to\mathbb{R}$ is defined by
$$
g(F,M_s):=\frac{\nu(1-\nu)\beta_E^2\nu_E^2F^2}{\alpha(F)((1-\nu)\nu_E\beta_EF+\alpha(F)\delta_M\gamma_sM_s)}
$$
(with $g(0,M_s):=0$),
{and the control $u:\mathbb{R}^+\to\mathbb{R}$ is a continuous function of time, where the notation~$\mathbb{R}^+$ stands for the interval $[0,+\infty).$
Since the biological model makes sense only for non-negative states,
we consider initial data $(F(0),M_s(0))\in\mathbb{R}^+\times\mathbb{R}^+$. Then,
being the function $g(F,M_s)$
globally Lipschitz-continuous on $\mathbb{R}^+\times\mathbb{R}^+$, we have existence and uniqueness of non-negative solutions for $t\in[0,T_{\max}]$, where
$$
T_{\max}:=\sup_{T\geq0}\{(F(t),M_s(t))\in\mathbb{R}^+\times\mathbb{R}^+\ \forall t\in[0,T]\}.
$$
Observe that positive invariance always holds for the $F$ component in \eqref{eq:reduced}.
Hence,  $T_{\max}$ is finite if  $M_s$ becomes negative, which can actually 
happen only if the control $u$ is negative somewhere.}
Now, under Assumption~\ref{ass:offspring}, the following stability result holds for~(\ref{eq:reduced}) in the absence of a control action (i.e.~$u\equiv0$), 
whose proof can be found in~\cite{almeida2022optimal}. 
\begin{prop}\label{eq:openloop}
Suppose that Assumption~\ref{ass:offspring} holds. Then, for $u(\cdot)\equiv0$, the reduced system \eqref{eq:reduced} admits two equilibria:
\begin{enumerate}
\item[i)] the extinction equilibrium $(F,M_s)=(0,0)$, which is unstable;\smallskip
\item[ii)] the persistence equilibrium $(F,M_s)=(\bar{F},0)$, with $\bar{F}=\frac{\nu\nu_Ek}{\delta_F}\left(1-\frac{1}{R_0}\right)$, which is locally asymptotically stable.
\end{enumerate}
\end{prop}
It is interesting to note that the instability of the extinction equilibrium cannot be proved by using the Lyapunov indirect method, i.e. by looking at  the linearized system, because the gradient of  $g(F,M_s)$ is {not continuous} at $(F,M_s)=(0,0)$ and therefore the linearized system is not well defined. {Nevertheless, the directional derivatives of $g(F,M_s)$ at $(F,M_s)=(0,0)$, restricted to directions in $\mathbb{R}^+\times\mathbb{R}^+$, are bounded.}

\section{Backstepping control}\label{sec:control}
Our goal is to design a (non-negative) control function $u(t)$ capable to stabilize the extinction equilibria. In particular we aim at design $u$ such that the closed-loop system is exponentially stable at the origin, in the sense of Definition~\ref{def:positive_stability} given next. The idea is to use a backstepping control design setup.
\begin{deff}[Positive\! Exponential\! Stabilization]\label{def:positive_stability}
{Given~the system~\eqref{eq:reduced} and a 
{locally Lipschitz-}continuous feedback control $u=\kappa(F,M_s)$, we say that $u$ globally positively exponentially stabilizes the system at $(F,M_s)=(0,0)$ with a rate $\lambda>0$ if there exists a constant $c_0>0$ such that, for any initial datum $(F(0),M_s(0))\in\mathbb{R}^+\times \mathbb{R}^+$,
one has {$T_{\max}=+\infty$ and}
\begin{equation}\label{expconv}
F^2(t)+M_s^2(t)\leq c_0 e^{-\lambda t}(F^2(0)+M_s^2(0))\ \ \forall t\geq0 
\end{equation}}
\end{deff}\smallskip
We first note that, due to $g(F,M_s)\geq0$, one has necessarily $F(t)\geq e^{-\delta_Ft}F(0)$ and, as such, there is no hope to obtain a decay rate for the female cluster larger than the inherent death rate $\delta_F$. 
{However, the next result  will entail that any decay rate $\delta_F-\epsilon>0$, with $\epsilon>0$ arbitrarily small, can be achieved
by choosing a suitable virtual control $M_s=M_s^*(F)$ for the first equation in \eqref{eq:reduced}. By the principle of backstepping, the goal will then be to design the actual control $u$ such that $M_s$ converges towards $M_s^*(F)$.}
\begin{prop}\label{prop:Mstar}
Let us pick $\hat{F}>\bar{F}$, where $\bar{F}$ corresponds to the persistence equilibrium for the female population, and define
\begin{equation}\label{eq:eps0}
{\epsilon:=\frac{\nu k\beta_E\nu_E}{\nu_E k+\hat{F}\beta_E+\delta_E k}}.
\end{equation}
Then the nonlinear feedback $M_s^*:\mathbb{R}^+\to\mathbb{R}$
defined by
\begin{equation}\label{eq:Mstar}
M_s^*(F):=\frac{(1-\nu)\nu_E\beta_E^2kF(\hat{F}-F)}{\gamma_s\delta_M(\beta_EF+k(\nu_E+\delta_E))^2}
\end{equation}
is positive for any $F\in(0,\hat{F})$, satisfies $M_s^*(0)=0$,
together with the equality $$g(F,M_s^*(F))=\epsilon F.$$
\end{prop}
\begin{proof}
{%
Let $\epsilon>0$.
By direct inspection, for $F>0$,
the equation $g(F,M_s)=\epsilon F$~rewrites as
$$
M_s=\frac{(1-\nu)\nu_E\beta_EkF(\frac{\nu k\beta_E\nu_E}\epsilon-(\beta_EF+\nu_Ek+\delta_E k))}{\gamma_s\delta_M(\beta_EF+k(\nu_E+\delta_E))^2}.
$$
For given $\hat{F}>0$, if $\epsilon$ is given by \eqref{eq:eps0} 
then the above right-hand side reduces to the function $M_s^*(F)$
defined in \eqref{eq:Mstar}, which vanishes at $F=0$ and
is positive 
for $F\in(0,\hat{F})$.
}
\end{proof}
\begin{remark}
{It is worth noticing that  the parameter $\epsilon$ tends to $0$ as $\hat{F}$ increases to $+\infty$,} thus enabling for the selection of a semiglobal
nonnegative virtual feedback law $M_2^*(F)$. 
\end{remark}

%
{From now on, let us fix $\hat{F}>\bar{F}$ sufficiently large such that the corresponding $\epsilon$ is smaller than $\delta_F$, and define the feedback law $M_s^*(F)$ accordingly.}\smallskip\\
Now, following the rationale of the backstepping approach \cite[Section 6.4]{isidori2017lectures}, the aim is to determine an actual control input $u$ capable of making the state $M_s$ track the desired feedback law $M^*_s(F)$ while, simultaneously, ensuring the overall system stability.
Towards this goal, let us first define the function $\pi(F, M_s):\mathbb{R}^+\times\mathbb{R}^+\to\mathbb{R}$ as 
\begin{equation}\label{eq:Pi}
\begin{array}{rl}
\pi(F, M_s)&:=\displaystyle\!\frac{g(F,M_s)\!-\!g(F,M^*_s(F))}{M_s-M^*_s(F)}\,F\smallskip\\
&\displaystyle=\frac{g(F,M_s)\!-\!\epsilon F}{M_s-M^*_s(F)}\,F\end{array}
\end{equation}
for $M_s\neq M^*_s(F)$, and 
$$\pi(F, M_s):=\frac{\partial g(F,M_s)}{\partial{M_s}}F$$
when $M_s=M^*_s(F)$.
{The function $\pi(F, M_s)$ 
\textcolor{black}{represents the mismatch rate between 
$g(F,M_s)F$ and $g(F,M_s^\star(F))F$, that we will need to absorb
through the control $u$; $\pi(F, M_s)$
is a well-defined and continuous function} on $\mathbb{R}^+\times\mathbb{R}^+$,
since $\frac{\partial g(F,M_s)}{\partial{M_s}}$ is bounded on 
$\mathbb{R}^+\times\mathbb{R}^+$ and it is everywhere continuous except at $(0,0)$. 


\begin{theo}\label{teo:back}
{Consider system \eqref{eq:reduced} with $u$ given by the feedback control
\begin{equation}\label{eq:backstepping}
\begin{array}{rl}
u^\star(F,M_s)&:=(\delta_s-\eta)M_s+\eta M^*_s(F)-\varrho\pi(F,M_s)\smallskip\\&+\displaystyle\frac{\partial M^*_s(F)}{\partial F}(g(F,M_s)-\delta_FF),
\end{array}
\end{equation}
for fixed gains $\varrho,\eta>0$ and \textcolor{black}{where $M^*_s(F)$ is given in \eqref{eq:Mstar} with $\epsilon>0$ defined by \eqref{eq:eps0} for 
given $\hat{F}>\bar{F}$.}
Then there exists $c_0>0$ such that 
any solution $(F,M_s)$ of \eqref{eq:reduced}
with $(F(0),M_s(0))\in\mathbb{R}^+\times\mathbb{R}^+$
satisfies \eqref{expconv} in its interval of existence,
with $\lambda=2\min\{\delta_F-\epsilon,\eta\}$.}
\end{theo}
\begin{proof}
Consider the Lyapunov function candidate
$$
\mathcal{V}(F,M_s):=\frac{\varrho}2 F^2+\frac12 (M_s-M_s^*(F))^2
$$
{which consists in the squared norm of $F$ and of virtual control mismatch $M_s-M^*_s(F)$, modulated by the parameter $\varrho$.}
Computing the derivative along the system trajectories, yields
$$
\begin{array}{rl}
\dot{\mathcal{V}}(F,M_s)&=\varrho Fg(F,M_s)-\varrho\delta_FF^2\smallskip\\
&+\displaystyle(M_s-M^*_s(F))\left(u-\delta_sM_s-\frac{\partial M^*_s(F)}{\partial F}\dot{F}\right)\smallskip\\
&=-\varrho(\delta_F-\epsilon)F^2+\varrho \pi(F,M_s)(M_s-M_s^*(F))\\
&+\displaystyle(M_s-M^*_s(F))\left(u-\delta_sM_s-\frac{\partial M^*_s(F)}{\partial F}\dot{F}\right).
\end{array}
$$
Now, replacing $u$ in the latter expression with the control law~\eqref{eq:backstepping}, one gets
$$
\begin{array}{rl}
\dot{\mathcal{V}}(F,M_s)&=-\varrho(\delta_F-\epsilon)F^2-\eta(M_s-M^*_s(F))^2\smallskip\\
&\leq -2\min\{\delta_F-\epsilon,\eta\}\mathcal{V}(F,M_s).
\end{array}
$$
{This implies that,
in the interval of existence of the solution, it holds that
$$\mathcal{V}(F(t),M_s(t))
\leq \mathcal{V}(F(0),M_s(0))e^{-\lambda t},$$ with 
$\lambda:=2\min\{\delta_F-\epsilon,\eta\}$.
From this, one deduces that the estimate \eqref{expconv}
holds with $c_0$ depending on $\rho$ only,
and with $(M_s-M^*_s(F))^2$ in place of $M_s^2$. 
Finally, using that the function $M_s^*(\cdot)$ satisfies $|M_s^*(F)|\leq K F$
for all $F\geq0$, for some $K>0$, together with some simple algebraic manipulations, 
one passes from the estimate on $(M_s-M^*_s(F))^2$ 
to the analogous one on $M_s^2$, with a $c_0$ also depending on $K$.
This concludes the proof.}
\end{proof}
{
Some comments about the controller~\eqref{eq:backstepping} are in order.
Firstly, we point out that it
provides exponential stabilization for the system, among solutions that
remains globally non-negative ($T_{\max}=+\infty$).
A sufficient condition for having $T_{\max}=+\infty$
is that $u\geq0$\footnote{This indeed implies the invariance of the set
$\mathbb{R}^+$ for $M_s$, and thus the global Lipschitz-continuity
of the function $g$.}.
However, the feedback control $u^\star$ defined by~\eqref{eq:backstepping}
may actually be negative somewhere,
{\color{black} which would be 
inconsistent with the biological model, as the control input in \eqref{eq:reduced} corresponds to the rate of release of sterile males and, as such, it should be non-negative}.
In the sequel, we will modify the definition of the
controller~\eqref{eq:backstepping} in order to cope with this problem.
This will lead to the global positive exponential stabilization of the system.
On the other hand, when $T_{\max}=+\infty$, 
the exponential convergence provided by Theorem~\ref{teo:back},
together with the Lipschitz-continuity of 
the controller $u^\star(\cdot,\cdot)$, which vanishes at $(0,0)$,
implies that $u^\star(F(t),M_s(t))$ also decays exponentially to $0$
as $t\to+\infty$. In particular, 
the total amount of released sterile mosquitoes is finite, as expressed by  $\int_0^{+\infty} u^\star(F(t),M_s(t))dt<+\infty$.\smallskip\\ } {Let us now investigate sufficient conditions for \eqref{eq:backstepping} to be indeed a non-negative feedback. We begin by observing that the function $\pi(F,M_s)$ is everywhere non-positive since the function $g(F,M_s)$ is non-increasing with respect to the second argument.} This observation allows us to neglect this term when evaluating the sign of $u^\star(F,M_s)$, in the sense that it is always helping the input to stay positive. Next result shows instead an interesting relationship between $M_s^*(F)$ and its derivative.
\begin{lemma}\label{lem:Mstar+}
Consider the function $M^*_s(F)$ defined in \eqref{eq:Mstar}. Then the following inequality holds
$$
M^*_s(F)\geq \frac{\partial M_s^*(F)}{\partial F}F\quad \forall F\in[0,\hat{F}]
$$
\end{lemma}
\begin{proof}
By direct inspection, one can check that
$$
\begin{array}{l}
M^*_s(F)- \displaystyle\frac{\partial M_s^*(F)}{\partial F}F\bigskip\\
\!\!=\displaystyle\frac{(1-\nu)\nu_E \beta_E^2kF^2(\beta_E(2\hat{F}-F)+k(\nu_E+\delta_E))}{\gamma\delta_M(\beta_E F+k(\nu_E+\delta_E))^3}
\end{array}
$$
which is clearly always non-negative in the interval of interest $F\in[0,\hat{F}]$.\end{proof}
\smallskip
Bearing the above properties in mind, we propose a corollary to Theorem~\ref{teo:back} that guarantees non-negativity of a slight variation of controller~\eqref{eq:backstepping}, which is still stabilizing whenever the gain $\eta$ is chosen in a suitable interval. \textcolor{black}{In particular, based on the previous arguments, it is clear that the only possibly negative term in the right-hand side of \eqref{eq:backstepping} is the last one. To rule out such an opportunity, before stating the result, it is useful to introduce the continuous function $(x,y)\in\mathbb{R}^2\mapsto \mathrm{cut}_2(x,y)\in\mathbb{R}$ with
$$
\mathrm{cut}_2(x,y):=\left\{
\begin{array}{rl}
0& \textrm{if}\  x<0\,\wedge\,y> 0\\
xy & \textrm{otherwise}
\end{array}
\right.
$$
that will be used to define the modified controller.
Such a function corresponds to the product $xy$ everywhere except for points lying in the $2^{\textrm{nd}}$ quadrant of the plane $(x,y)$, where it is identically zero.}
\begin{corol}\label{cor:positivity}
{The conclusion of Theorem \ref{teo:back}
holds true with the controller \eqref{eq:backstepping}
replaced by
}
\begin{equation}\label{eq:upiu}
\begin{array}{rl}
u^\star_+(F,M_s)&:=(\delta_s-\eta)M_s+\eta M^*_s(F)-\varrho{\pi(F,M_s)}\smallskip\\&+\displaystyle\mathrm{cut}_2\left(\frac{\partial M^*_s(F)}{\partial F},g(F,M_s)-\delta_FF\right)
\end{array}
\end{equation}
{Moreover, if $\eta\in(\delta_F,\delta_s)$, 
then $u^\star_+(F,M_s)\geq 0$ for any $(F,M_s)\in[0,\hat{F}]\times\mathbb{R}^+$.}
\end{corol}
\begin{proof}
{Let us first check that, under the assumption
{\color{black}$\eta\in(\delta_F,\delta_s)$}, the control $u^\star_+(F,M_s)$ is non-negative for any $(F,M_s)\in[0,\hat{F}]\times\mathbb{R}^+$. Thanks to the condition $\pi(F,M_s)~\leq~0$,
the only case where such a property may fail
is when the last term in \eqref{eq:upiu} is negative, that is, by
the definition of $\mathrm{cut}_2$,
when $$\frac{\partial M^*_s(F)}{\partial F}>0\ \wedge\ g(F,M_s)-\delta_FF<0.$$
Observe that the last term is equal to $\dot{F}$.
On one hand, recalling the definition of $M^*_s(F)$, there exists $F_1\in(0,\hat{F})$ such that $\frac{\partial M^*_s(F)}{\partial F}$ is positive only for $F\in(0,F_1)$. On the other hand, the bound $\dot{F}\geq -\delta_F F$
holds for any $(F,M_s)\in\mathbb{R}^+\times\mathbb{R}^+$ and then, invoking Lemma~\ref{lem:Mstar+}, we find 
$$-\frac{\partial M^*_s(F)}{\partial F}\dot{F}\leq
\delta_F M^*_s(F)$$ for all $(F,M_s)\in[0,\hat{F}]\times\mathbb{R}^+$.
Summing up, using $\eta>\delta_F$, we derive}
$$
\eta M^*_s(F)+\mathrm{cut}_2\left(\frac{\partial M^*_s(F)}{\partial F},\dot{F}\right)\geq0\quad \forall (F,M_s)\in\mathbb{R}^+\times\mathbb{R}^+ 
$$
thus showing non-negativity of $u_+^\star$.\smallskip\\ To prove stabilization, 
consider again the Lyapunov function candidate used in the proof of Theorem~\ref{teo:back} and evaluate its derivative along the solution of~\eqref{eq:reduced} driven by the control $u^\star_+$. This yields
$$
\begin{array}{rl}
\dot{\mathcal{V}}(F,M_s)&\leq-\varrho(\delta_F-\epsilon)F^2-\eta(M_s-M^*_s(F))^2\vspace{0.2cm}\\
&\!\!\!\!\!\!\!\!\!\!\!\!\!\!\!\!\!\!\!+\displaystyle(M_s-M^*_s(F))\left(\!\mathrm{cut}_2\left(\frac{\partial M^*_s(F)}{\partial F},\dot{F}\right)
-\frac{\partial M^*_s(F)}{\partial F}\dot{F}\!\right)
\end{array}
$$
Now, when either $\frac{\partial M^*_s(F)}{\partial F}\dot{F}\geq0$ or $\dot{F}<0$, the last term in the right-hand side is zero and we get back to the original condition obtained in the proof of Theorem~\ref{teo:back}. Conversely, when $$\frac{\partial M^*_s(F)}{\partial F}\dot{F}<0 \ \wedge\ \dot{F}>0,$$
the additional term $$(M_s^*(F)-M_s)\frac{\partial M^*_s(F)}{\partial F}\dot{F}$$ enters the Lyapunov inequality. 
It is straightforward to see that, whenever $M_s^*(F)\geq M_s$, the latter term is negative and thus may indeed be helpful for the convergence. 
Furthermore, we can actually prove that this is always the case. Indeed, for $M_s^*(F)< M_s$, one would have $g(F,M_s)<\epsilon F$ and so $\dot{F}<(\epsilon-\delta_F)F\leq 0$, thus contradicting the hypothesis $\dot{F}>0$.
In conclusion we have shown that, in all possible scenarios, the Lyapunov inequality
$$\begin{array}{rl}
\dot{\mathcal{V}}(F,M_s) &\leq-\varrho(\delta_F-\epsilon)F^2-\eta(M_s-M^*_s(F))^2\\
&\leq-2\min\{\delta_F-\epsilon,\eta\}\mathcal{V}(F,M_s)
\end{array}
$$
holds true, thus proving that the equilibrium $(0,0)$ is still exponentially stable, with a decay rate not smaller than $2\min\{\delta_F-\epsilon,\eta\}.$
\end{proof}
\begin{remark}\label{rem:noneed}
It is worth stressing that the condition introduced in Corollary~\ref{cor:positivity} might be slightly conservative. This will be illustrated in the simulation results, where the ``pure'' backstepping control $u^\star$ is naturally non-negative without the need of introducing the $\mathrm{cut}_2$ correction.
\end{remark}
{The feedback controller $u_+^\star$ defined
in Corollary \ref{cor:positivity}, with 
$\rho>0$ and $\eta\in(\delta_F,\delta_s)$,
is non-negative and well defined for $(F,M_s)\in[0,\hat{F}]\times\mathbb{R}^+$.\smallskip\\ To get a result holding globally for $F\in\mathbb{R}^+$, we further modify $u_+^\star$ as follows.
Fix a number $F_2\in(\bar F,\hat{F})$~\footnote{One could take for instance
$F_2=(\bar F+\hat{F})/2$.}
and consider a cut-off
function $\chi(F)$ which is smooth, non-increasing
and satisfies $\chi(F)\!=\!1$ for $F\!\leq\! F_2$ and
$\chi(F)\!=\!0$ for $F\!\geq\! \hat{F}$.
We then~set 
\begin{equation}\label{eq:utilde}
\tilde u^\star(F,M_s):=u_+^\star(F,M_s)\chi(F)\end{equation}
for $(F,M_s)\in\mathbb{R}^+\times\mathbb{R}^+$.
This control is non-negative.
This implies in particular
that, under its action, the set $\mathbb{R}^+$ is invariant for the component
$M_s$ of the solution to \eqref{eq:reduced}, i.e. $T_{\max}=+\infty$ (recall that  $F$ always remains non-negative).\smallskip\\ 
Next, in order to evaluate the system behaviour subject to
the control $\tilde u^\star$, let us analyze the evolution of $F$ in the region $F\geq F_2$.}
Using the monotonicity of $g(F,M_s)$ in the second variable,
we find
$$g(F,M_s)-\delta_F F\leq g(F,0)-\delta_F F
=\left(\frac{\nu\beta_E\nu_E}{\alpha(F)}-\delta_F\right)F$$
for all $(F,M_s)\in\mathbb{R}^+\times\mathbb{R}^+$.
Then, since $\alpha(F)$ is increasing, we have
$g(F,M_s)-\delta_F F\leq-\sigma F$ for $F\geq F_2$, with 
\begin{equation}\label{eq:sigma}
\sigma:=\delta_F-\frac{\nu\beta_E\nu_E}{\alpha(F_2)}\end{equation}
and one can check that $F_2>\bar F$ implies $\sigma>0$.
As a consequence, for a solution $(F,M_s)$ of
\eqref{eq:reduced}, it holds that 
$\dot F\leq -\sigma F$ whenever $F\geq F_2$.
This shows, on one hand, that $[0,F_2]$ is an invariant set for the $F$ component,
and on the other hand that if $F(0)>F_2$ then $F(t)$ decreases exponentially, 
with a rate not smaller than $\sigma$,
{until a time $T$ where $F(T)=F_2$.
After the transient time $T$  ($=0$ if $F(0)\leq F_2$), 
the expression of $\tilde{u}^\star(F,M_s)$ boils down 
to~\eqref{eq:upiu}.
It then follows from Corollary~\ref{cor:positivity} that, for $t\geq0$,
%
\begin{equation*}
F^2(t)\leq c_0 e^{-\lambda' t}(F^2(0)+M_s^2(0)) 
\end{equation*}
with $\lambda':=\min\{2\sigma,\lambda\}$.
Next, using the fact that the controller $\tilde{u}^\star$ satisfies
$\tilde{u}^\star(M_s,F)\leq(\delta_s-\eta) M_s+KF$ for all $(F,M_s)\in\mathbb{R}^+\times\mathbb{R}^+$, for some positive constant $K$,
one has 
$$\dot{M_s}=\tilde{u}^\star-\delta_s M_s\leq 
-\eta M_s+KF.$$
From this, using the exponential bound for $F$ obtained before, 
one infers
$
M_s^2(t)\leq c_1 e^{-\lambda^{\prime\prime} t}(F^2(0)+M_s^2(0)) 
$
with $\lambda^{\prime\prime}:=\min\{2\eta,\lambda'\}$
and $c_1$ depending on $c_0,K,\eta,\lambda$.}\smallskip\\
In conclusion, summarizing the above discussion, we can state the following global stabilization result.
\begin{theo}
Let us consider the reduced mosquito population dynamics \eqref{eq:reduced}. The non-negative feedback $\tilde{u}^\star(F,M_s)$ defined by \eqref{eq:backstepping}-\eqref{eq:upiu}-\eqref{eq:utilde},
{%
with $\rho>0$ and $\eta\in(\delta_F,\delta_s)$}, guarantees the global positive exponential stabilization of the extinction equilibrium $(F,M_s)=(0,0)$, with convergence rate not smaller than $\lambda=2\min\{\delta_F-\epsilon,\eta,\sigma\}$ 
where $\sigma<\delta_F$ is defined in \eqref{eq:sigma}.   
\end{theo}

\section{Simulations}\label{sec:simul}
The efficacy of the backstepping control strategy \eqref{eq:Mstar}-\eqref{eq:backstepping} has been successfully tested and validated in simulations {on the reduced model \eqref{eq:reduced} and, without any formal stability or performance guarantee yet, also on the complete model \eqref{eq:fullmodel}. As a comparison, we have also simulated the response of the complete model~\eqref{eq:fullmodel} to the feedback control obtained in \cite[Eq. (7) and Remark 3.4]{almeida2022optimal}.}\smallskip\\ 
We have considered the following values for the mosquitoes population parameters, \textcolor{black}{consistently with \cite{strugarek2019use}:}
\begin{table}[htp]
\begin{center}
\begin{tabular}
{|c|c|c|c|c|c|c|c|}
\hline
$\beta_E$ & $\gamma_s $&$\nu_E$&$\nu$&$\delta_E$&$\delta_M$&$\delta_F$&$\delta_s$\\
\hline
10&1&0.005&0.49&0.03&0.1&0.04&0.12\\
\hline
\end{tabular}
\end{center}
\caption{Nominal model parameters}\label{tab:values}
\end{table}%

and the following control parameters:
$$
\epsilon=0.01,\ \hat{F}\simeq\frac{27}{20}\bar{F},\ \eta=\delta_s-0.02,\ \varrho=0.5
$$
All simulations have been performed initializing the system at the persistence equilibrium $(\bar{E},\bar{M},\bar{F},0)$ for the complete model
with
$$
\bar{M}=5106,\ \bar{E}=200240,\ \bar{F}=12264,$$
and \textcolor{black}{environmental capacity $k=\left(1-\frac{\delta_F(\nu_E+\delta_E)}{\beta_E\nu\nu_E}\right)^{-1}
\!\bar{E}=212370$.}\smallskip\\ 
The response of the system to the backstepping control policy is illustrated in the next figures. {In particular, Figure~\ref{fig:female} shows the evolution of the females: the behaviour of this compartment of mosquitoes population, driven by the control \eqref{eq:Mstar}-\eqref{eq:backstepping}, is quite similar for the reduced model and the complete model, in both cases being characterized by a clear exponential decay. In contrast to such similarity, the behavior of  control inputs shows a remarkable difference as visible in Figure~\ref{fig:control}. In particular, in the case of complete dynamics, the initial behaviour is equivalent to the one of reduced dynamics but then the decay becomes slower, thus indicating that a larger number of sterile males is needed to guarantee the reduction of females when their number decreases.} It is also interesting to notice that the control input is vanishing and always non-negative (see Remark~\ref{rem:noneed}), thus showing the efficacy of the proposed backstepping approach in driving the female mosquitoes towards the extinction equilibrium.\smallskip\\  The obtained results suggest that the feedback control law \eqref{eq:Mstar}-\eqref{eq:backstepping}, although being tailored for the reduced model, is actually applicable with success also for the stabilization of the complete model.
{Conversely, it is apparent that the feedback law in \cite[Eq. 7]{almeida2022optimal}, which was specifically designed for the reduced model with optimality purposes, is instead not suitable for the stabilization of the complete system.
In fact, on the one hand, persistency of oscillations is visible in the dashed-dotted plots in Figure~\ref{fig:female}-\ref{fig:males} and, on the other hand, a non-vanishing control behaviour characterizes the dashed-dotted plots in Figure~\ref{fig:control}.} \smallskip\\
{To better illustrate the capabilities of the proposed control law, we have also performed a robustness check. In particular the nominal control law, defined based on the parameter values given in Table~\ref{tab:values}, has been fed to a system whose parameters were subject to a 10\% random uncertainty with respect to the nominal values in Table~\ref{tab:values}. The results obtained for $20$ random iterations, reported in Figures~\ref{fig:female_rob}-\ref{fig:control_rob}, are quite promising and support the claim of robustness, with stabilization of the extinction equilibrium for the uncertain system even when the latter is controlled by the nominal feedback law. It is worth noticing that the uncertainty in the model may lead to larger values of the control function compared to the nominal case (see Figure~\ref{fig:control_rob}). Nevertheless, we can observe that the asymptotically decreasing behaviour of the control is preserved in all the performed iterations. It is also foreseen that better performances are likely to be obtained when the control parameters are selected according to a worst-case choice given a certain bounded range of model uncertainties, thus paving the road for a future formal proof of robust exponential stabilization by means of the proposed backstepping controller.}

\begin{figure}[b!]
\centering
\includegraphics[width=0.88\columnwidth]
{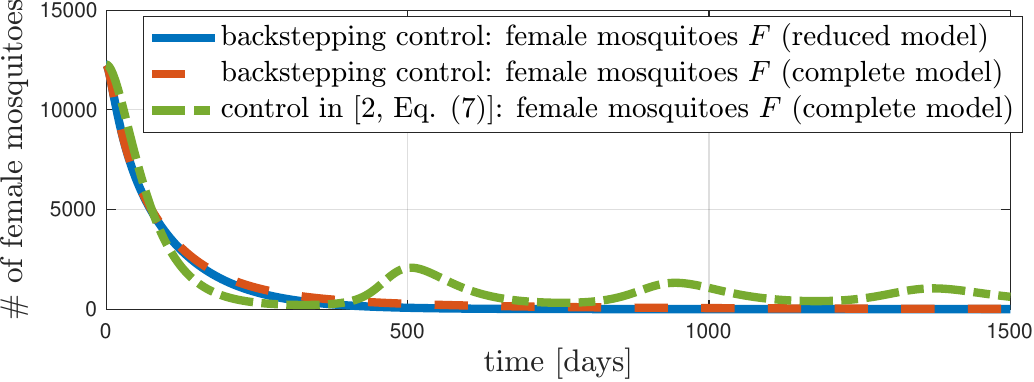}
\caption{Female mosquitoes under backstepping control}\label{fig:female}
\end{figure}

\begin{figure}[b!]
\centering
\includegraphics[width=0.88\columnwidth]
{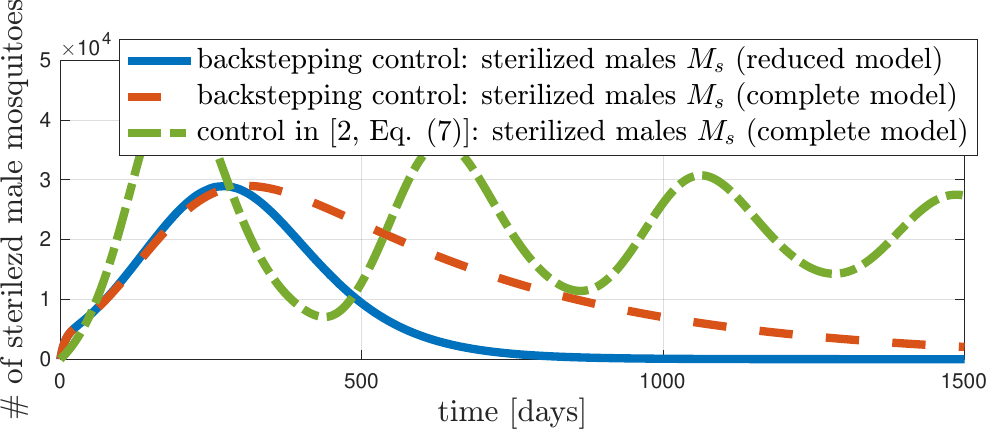}
\caption{Sterilized male mosquitoes under backstepping control}\label{fig:males}
\end{figure}

\begin{figure}[b!]
\centering
\includegraphics[width=0.88\columnwidth]
{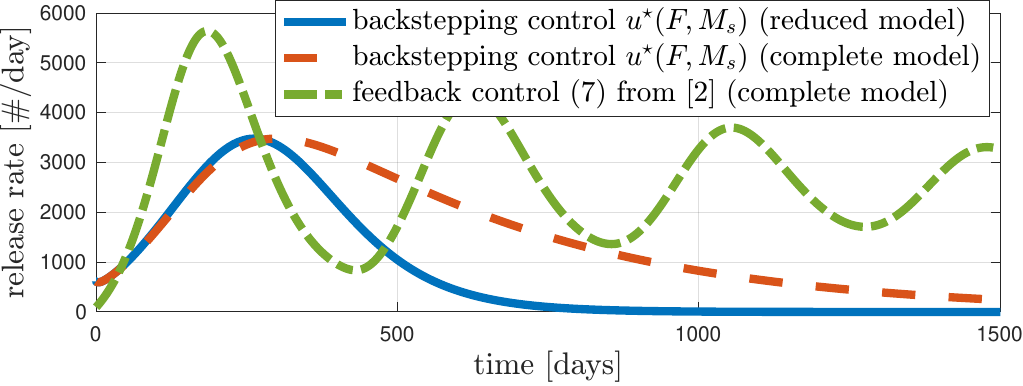}
\caption{Feedback control functions}\label{fig:control}
\end{figure}

\section{Conclusions and discussion}\label{sec:conclusion}
In this paper we have considered a stabilization problem of the extinction equilibrium for a population of mosquitoes through the sterile insect technique. Using a model-based setup, a backstepping control law has been considered and proved to guarantee the global stabilization of a reduced model comprising females and sterile males only. Sufficient conditions are provided for the non-negativity of the control function, representing the rate of release of sterile males. Furthermore, the same feedback control law has been tested in simulations for the complete model, showing its efficacy for this case too. {A preliminary simulation study on the inherent robustness of the proposed backstepping control has been also performed with success, both for the reduced and the complete model.} Motivated by the encouraging testing in simulations, currently we are working on enhancing the proposed backstepping control approach by formally considering robustness to {uncertainties in the model parameters \textcolor{black}{and in the environmental capacity~$k$} as previously addressed in, {\it e.g.}, \cite{bliman2022robust},  and introducing sampled measurements, to better catch the typical discontinuities in data acquisition and during the releasing process}~\cite{bidi2024feedback}. Moreover, we aim at formally proving global stability for the complete model under the proposed backstepping controller. Future developments will also be devoted to dealing with diffusion terms in the model, similarly to what has been done in \cite{almeida2022wave},~\cite{almeida2023analysis},\cite{leculier2023control}, and {to consider control based on the incompatible insect technique~\cite{zheng2019incompatible},~\cite{bliman2019feedback}.}
\begin{figure}[b!]
\centering
\includegraphics[width=0.88\columnwidth]
{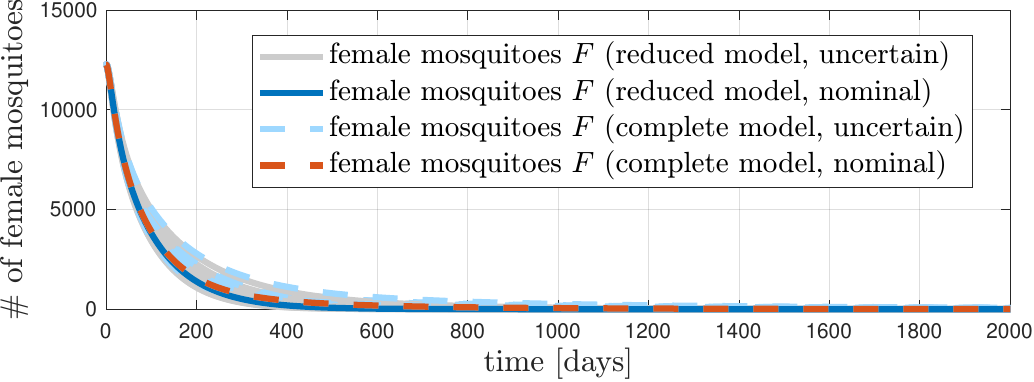}
\caption{Female mosquitoes under backstepping control, robustness analysis: nominal system VS uncertain system (20 iterations)}\label{fig:female_rob}
\end{figure}

\begin{figure}[b!]
\centering
\includegraphics[width=0.88\columnwidth]
{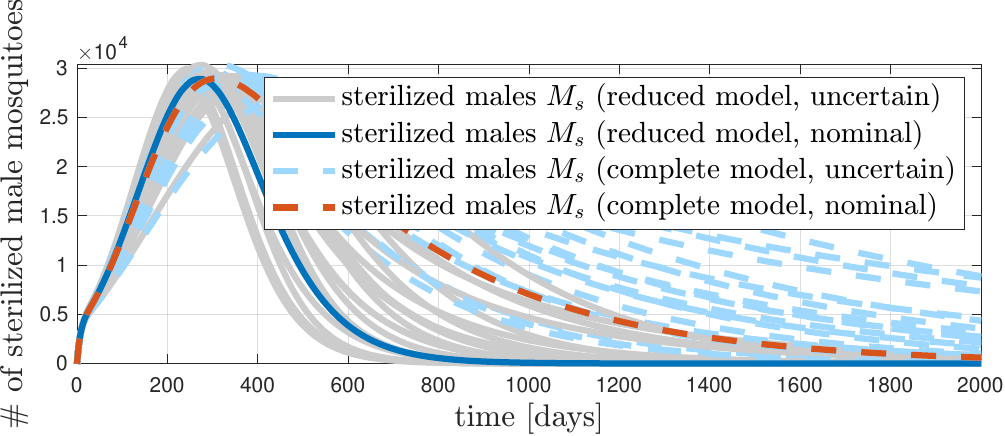}
\caption{Sterilized male mosquitoes under backstepping control, robustness analysis: nominal system VS uncertain system (20 iterations)}\label{fig:males_rob}
\end{figure}

\begin{figure}[b!]
\centering
\includegraphics[width=0.88\columnwidth]
{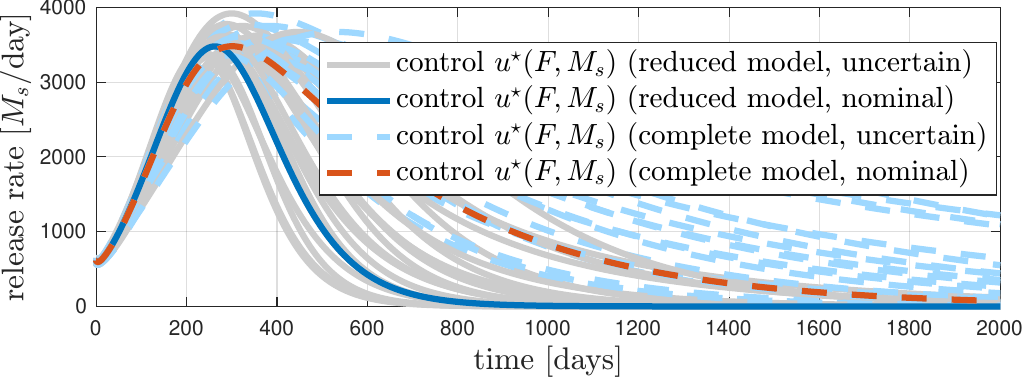}
\caption{Feedback control functions, robustness analysis: nominal system VS uncertain system (20 iterations)}\label{fig:control_rob}
\end{figure}




\bibliography{mosquito}

\end{document}